# A GENERALIZATION OF DESARGUES' INVOLUTION THEOREM

NICHOLAS PHAT NGUYEN

**Abstract**. This paper states and proves a generalization of the well-known Desargues' involution theorem from plane projective geometry.[1]

**1. INTRODUCTION.** A classic theorem in projective geometry is the following remarkable result discovered by Girard Desargues, a French engineer and mathematician widely regarded as one of the founders of projective geometry.

**Theorem (Desargues' Involution Theorem):** *Consider four points in general position in the real projective plane, i.e., no three of these four points are collinear. Let $\mathcal{F}$ be the family of conics passing through these four points. Then for any line $\ell$ that does not pass through any of these four points, each conic in $\mathcal{F}$ will, if it intersects $\ell$, do so in a pair of points that are conjugate under an involution of the line $\ell$.*

Here, an involution means a projective transformation of the line $\ell$ that has order two. A conic is the set of zero points of a quadratic form $\neq 0$ in projective dimension two. Specifically, if the given projective plane represents the set of lines in a vector space $V$ of linear dimension 3, then a conic is the set of isotropic points of a symmetric bilinear form $q$ defined on $V$ i.e., the projective set associated with vectors $v$ such that $q(v, v) = 0$. Such a set of isotropic or zero points of a quadratic form is generally called a quadric, so a conic is the special case of a quadric in projective dimension two.

---

In this note, we want to prove the following generalization of the Desargues' involution theorem.

**Theorem (Generalization of Desargues' Involution Theorem):** *Consider a projective space of any dimension over a field K of characteristic ≠2. A pencil of quadrics in that projective space will intersect a line in that space (when they do) in pairs of points that are conjugate under an involution if and only if (i) the line is not part of a quadric in that pencil and (ii) the line contains no common zero point for all the quadrics in that pencil. Equivalently, these conditions are satisfied if and only if the pencil contains two quadrics whose intersection points with the given line over an algebraic closure of K are disjoint.*

The term "pencil" in this context means a linear system of projective dimension 1. In other words, if $R$ and $S$ are two symmetric bilinear forms that represent two different quadrics in that pencil, then all of the quadrics in the pencil can be represented by the symmetric bilinear forms $aR + bS$, where $a$ and $b$ are numbers running through the ground field $K$.

As this generalization suggests, the Desargues' involution theorem is essentially an algebraic property of the projective line and has little to do with the ambient projective space, and does not require any continuity to hold true.

The famous butterfly theorem of Euclidean plane geometry is a special case of the Desargues' involution theorem. With our generalization of the Desargues' involution theorem, we will also state a more general form of the butterfly theorem for affine space of any dimension over a field of characteristic $\neq 2$.

**2. SYMMETRIC BILINEAR FORMS.** We will work with a ground field $K$ of characteristic $\neq 2$. A projective line defined over $K$ is the set $P(E)$ of lines through the origin in a $K$-vector space $E$ of linear dimension 2. If we identify such a vector space with the affine plane $K^2$ through a suitable choice of basis, we can describe the projective line as $K \cup \{\infty\}$, i.e., the affine line $K$ extended by adjoining a point at infinity denoted by the symbol $\infty$, where each element $x$ of $K$ is identified with the line in $K^2$ passing through the point $(x, 1)$ and the origin, and the point ∞ is identified with the horizontal line consisting of all the points ($u$, 0) (with $u$ running through $K$) in $K^2$.

A symmetric bilinear form on $E$ is a mapping from the product space $ExE$ into the ground field $K$ that is linear in each variable and symmetric. Relative to a given basis of $E$, such a pairing is described by a symmetric 2 by 2 matrix, and so the set $B$ of all symmetric bilinear forms naturally has the structure of a $K$-vector space of linear dimension 3. A symmetric bilinear form on $E$ is non-degenerate if and only if the determinant of the corresponding matrix is non-zero.

Given a symmetric 2 by 2 matrix $\begin{pmatrix} a & b \\ b & c \end{pmatrix}$, its determinant $ac - b^2$ is a quadratic form of the matrix entries, so we can associate a symmetric bilinear form to the determinant. Specifically, we can define the following pairing $BxB \to K$ for any two symmetric 2 by 2 matrices:

$$<\begin{pmatrix} a & b \\ b & c \end{pmatrix}, \begin{pmatrix} a' & b' \\ b' & c' \end{pmatrix}> = \tfrac{1}{2}(ac' + ca') - bb'$$

Such a pairing is well-defined because the ground field $K$ has characteristic ≠ 2, and it is clearly bilinear and symmetric. With this pairing, $B$ is a quadratic space isometric to the orthogonal sum of an Artinian plane (also known as a hyperbolic plane) and the one-dimensional space $K$ (with multiplication). Accordingly, this pairing on $B$ is regular, meaning that $<x, y> = 0$ for all $y$ in $B$ if and only if $x = 0$. The non-isotropic elements of $B$ are the symmetric matrices $z$ such that $<z, z> = \det z \neq 0$, i.e., the matrices that represent non-degenerate symmetric bilinear forms on $E$.

Because the determinant pairing is regular, any non-degenerate symmetric bilinear form $q$ gives us an orthogonal decomposition $B = Kq \oplus$ (the orthogonal complement of $Kq$). The orthogonal complement of $Kq$ has linear dimension 2, so whenever there is a subspace $W$ of $B$ of dimension at least 2 such that $W$ is orthogonal to a non-degenerate element $q$ of $B$, then $W$ must in fact be the 2-dimensional orthogonal complement of $Kq$, and $B$ must be the direct sum of $Kq$ and $W$.

Implicit in the determinant pairing is the choice of a basis for $E$ to allow us to identify the space $B$ with the space of 2 by 2 symmetric matrices. If we change the coordinates for $E$ by a general linear transformation $x = S(x')$, then the matrix $M$ of a symmetric bilinear form in the old coordinate $x$ will become ${}^tSMS$ (where ${}^tS$ is the transpose of $S$) in the new coordinates $x'$. If we take the norm $<M, M>$ of $M$ under the determinant pairing relative to the old basis, we get $\det(M)$. On the other hand, the norm of $M$ in the determinant pairing relative to the new basis becomes $<{}^tSMS, {}^tSMS> = \det({}^tSMS) = det(S)^2\det(M)$. Therefore the simple linear transformation $q \mapsto (\det S)q$ gives us an isometry between the vector space $B$ with the determinant pairing in the new coordinate $x'$ and the vector space $B$ with the determinant pairing in the old coordinate $x$. Accordingly, orthogonal properties in the space $B$ under the determinant pairing are independent of any basis chosen for the coordinates.

A symmetric bilinear form $q \neq 0$ on $E$ falls into one of 3 types:

- The form $q$ is non-degenerate and anisotropic, meaning there is no non-zero vector $v$ in $E$ such that $q(v, v) = 0$. The form $q$ has no zero or isotropic point on the projective line $P(E)$.
- The form $q$ is non-degenerate and isotropic. In that case, the quadratic space $(E, q)$ is isometric to an Artinian plane, and the isotropic vectors $v$ such that $q(v, v) = 0$ consist of exactly two different lines in $E$. The form $q$ has two distinct isotropic points on the projective line $P(E)$. Any symmetric bilinear form that has the same isotropic vectors as $q$ must be a scalar multiple of $q$. Indeed, relative to the basis consisting of those two isotropic vectors the matrices of these two bilinear forms both have zeros in the diagonal and a non-zero number in the cross diagonal.
- The form $q$ is degenerate. In that case, its radical is a one-dimensional subspace of $E$, and the form q has exactly one isotropic point on the projective line $P(E)$. Any symmetric bilinear form that has the same isotropic vectors as $q$ must be a scalar multiple of $q$. Indeed, relative to any basis that includes an isotropic vector, the matrices of these two bilinear forms both have just one non-zero number in the same diagonal position.

The above three situations tell us how a projective line in an arbitrary projective space can intersect a quadric. Let such a quadric be the projective set of isotropic points of a symmetric bilinear form $q \neq 0$ defined on a vector space $V$, and let the given projective line be $P(E)$, where $E$ is a subspace of dimension 2 of $V$. If the restriction of $q$ to $E$ is a non-zero symmetric bilinear form on $E$, then: (i) such a restriction being anisotropic means that the line $P(E)$ does not intersect the given quadric; (ii) such a restriction being non-degenerate and isotropic means that the line $P(E)$ intersects the given quadric in two distinct points; and (iii) such a restriction being degenerate means that the line $P(E)$ intersects the given quadric in one tangent point. Of course, if the restriction of $q$ to $E$ is the zero symmetric bilinear form on $E$, then that simply means the entire line $P(E)$ is part of the quadric.

When the form $q$ has one or two isotropic points we can tell what the symmetric matrices that represent $q$ (relative to a given basis of $E$) would look like. Suppose that relative to a given basis of $E$, the form $q$ has isotropic vectors $(s, t)$ and $(u, v)$ which are linearly independent in the non-degenerate case and proportional in the degenerate case. The following quadratic form has the same isotropic vectors, and therefore must be the same quadratic form associated with $q$, up to a scalar factor:

$$(tX - sY)(vX - uY) = tvX^2 - (tu + sv)XY + suY^2$$

Accordingly, the 2 by 2 symmetric matrix associated with $q$ would look like the following, up to a scalar factor:

$$\begin{pmatrix} tv & w \\ w & su \end{pmatrix}, \text{ where } w = -\tfrac{1}{2}(tu + sv).$$

3. **INVOLUTIONS.** A projective transformation of a projective line $P(E)$ is induced by an invertible linear transformation of $E$. So relative to a given basis of $E$, a projective transformation can be described by an invertible 2 by 2 matrix with coefficients in the ground field $K$, up to a scalar factor. A major focus of the classical geometry of projective line is the study of involutions, defined as projective transformations of order two. For a projective transformation to be an involution, it is necessary and sufficient that any 2 by 2 matrix $T$ corresponding to such an involution satisfies a minimal equation $T^2 - s = 0$ with $s$ $\neq 0$. At the same time, we know from the Hamilton-Cayley theorem that $T^2 - \text{tr}(T)T + \det(T) = 0$. So involutions are given by invertible 2 by 2 matrices with zero trace.

Given an involution on $P(E)$ represented by $\begin{pmatrix} a & b \\ c & -a \end{pmatrix}$, we can map that to the element of $P(B)$ represented by $\begin{pmatrix} -c & a \\ a & b \end{pmatrix}$. Note that the determinant of $\begin{pmatrix} -c & a \\ a & b \end{pmatrix}$ is equal the determinant of $\begin{pmatrix} a & b \\ c & -a \end{pmatrix}$, so that the symmetric matrix $\begin{pmatrix} -c & a \\ a & b \end{pmatrix}$ represents a non-degenerate symmetric bilinear form, or a non-isotropic element in the quadratic space $B$ with the determinant pairing.

This correspondence is a bijection between involutions on $P(E)$ and non-degenerate elements of $P(B)$, with the inverse correspondence mapping a non-degenerate element of $P(B)$ represented by $\begin{pmatrix} a & b \\ b & c \end{pmatrix}$ to the involution on $P(E)$ represented by $\begin{pmatrix} b & c \\ -a & -b \end{pmatrix}$.

We will refer to the element of $P(B)$ represented by the matrix $\begin{pmatrix} -c & a \\ a & b \end{pmatrix}$ as the Desargues bilinear form associated with the involution represented by $\begin{pmatrix} a & b \\ c & -a \end{pmatrix}$. While this bijective correspondence between involutions and non-degenerate elements of $P(B)$ is defined based on the matrix representation of involutions and symmetric bilinear forms relative to a specific basis of the vector space $E$, the correspondence is really independent of any basis, based on the following proposition.

**Proposition 1:** *Each pair of conjugate points in an involution are the isotropic points of a symmetric bilinear form that is orthogonal (under the determinant pairing) to the Desargues bilinear form corresponding to the given involution.*

*Proof.* Consider an involution on $P(E)$ represented by $\begin{pmatrix} a & b \\ c & -a \end{pmatrix}$, and let $(x, y)$ and $(ax + by, cx - ay)$ be a generic pair of conjugate points. These points can be proportional if $(x, y)$ happens to be a fixed point of the involution. Note that if we take any two conjugate points of the given involution and treat them as isotropic points of a symmetric bilinear form, then that form is uniquely determined up to a scalar factor, per our discussion of symmetric bilinear forms above. A symmetric bilinear form with these two isotropic points has the form:

$$\begin{pmatrix} y(cx - ay) & z \\ z & x(ax + by) \end{pmatrix}, \text{ where } z = -\tfrac{1}{2} y(ax + by) - \tfrac{1}{2} x(cx - ay).$$

The determinant pairing between that bilinear form and the Desargues form $\begin{pmatrix} -c & a \\ a & b \end{pmatrix}$ gives us the number:

$$\tfrac{1}{2}(-cx(ax + by) + by(cx - ay)) - az =$$

$$= \tfrac{1}{2}(-cx(ax + by) + by(cx - ay)) + \tfrac{1}{2}(ay(ax + by) + ax(cx - ay)) =$$

$$= \tfrac{1}{2}(-caxx - bayy) + \tfrac{1}{2}(abyy + acxx) = 0. \blacksquare$$

Recall that orthogonal relationships under the determinant pairing is independent of any basis. Hence, regardless of any basis we choose for $E$, the Desargues bilinear form corresponding to a given involution is orthogonal to all the elements of $P(B)$ associated with pairs of conjugate points of the involution.

Because we can find at least two different pairs of such conjugate points, these elements of $P(B)$ generate a subspace of projective dimension at least one. Since a Desargues bilinear form is non-degenerate by construction, the above proposition implies that all the symmetric bilinear forms associated with pairs of conjugate points of the involution generate a projective subspace of dimension 1 (a pencil) in $P(B)$, and the Desargues bilinear form is the uniquely determined bilinear form up to a scalar factor that is orthogonal to that pencil.

**4.    GENERALIZATION OF DESARGUES' INVOLUTION THEOREM.** Based on Proposition 1, there is a natural bijection (independent of any basis or coordinates) between involutions of a projective line and non-degenerate bilinear forms on that line or, equivalently, the orthogonal complements of such non-degenerate bilinear forms. A necessary and sufficient condition for a pencil of bilinear forms on a projective line to be the orthogonal complement of a non-degenerate bilinear form is that the pencil must be regular under the determinant pairing. In other words, that pencil must correspond to a two-dimensional subspace of $B$ which is a regular quadratic space under the determinant pairing.

In general, it is straight-forward to check whether a bilinear space of dimension 2 is regular under a given symmetric pairing. We can just write down the matrix of that pairing relative to a suitable basis of the space and determine if the matrix has non-zero determinant. In our particular case, we also have another geometric criterion.

**Proposition 2:** *A 2-dimensional subspace of the space of all symmetric bilinear forms on a vector space of dimension 2 is regular with respect to the determinant pairing if and only if there is no common isotropic vector for all the forms in that subspace, or equivalently, if the subspace can be generated by two forms with no common isotropic vector.*

*Proof.* As before, we denote by $B$ the space of all symmetric bilinear forms on a 2-dimensional $K$-vector space $E$. With a choice of basis for $E$, $B$ can be identified with the space of symmetric 2 by 2 matrices with coefficients in $K$ and therefore has linear dimension 3. Let $G$ be a 2-dimensional subspace of $B$. Recall that $B$ is regular under the determinant pairing. Because $B$ is regular, the subspace $H$ of $B$ orthogonal to $G$ is therefore a 1-dimensional subspace, say generated by a bilinear form $h$.

$G$ is a regular subspace of $B$ if and only if $h$ does not belong to $G$, i.e., if and only if $h$ is non-isotropic. If $h$ is isotropic, then it is a degenerate bilinear form orthogonal to all the forms in in $G$. Moreover, $h$ has a one-dimensional radical. Choose a basis of $E$ where the first vector of that basis lies in the radical of $h$. Relative to that basis, the matrix of $h$ has the form $\begin{pmatrix} 0 & 0 \\ 0 & x \end{pmatrix}$ where $x \neq 0$. Any symmetric bilinear form $\begin{pmatrix} a & b \\ b & c \end{pmatrix}$ orthogonal to $h$ under the determinant pairing must satisfy the equation $ax = 0$, i.e., $a = 0$. That means the first vector of the basis is also an isotropic vector for any form orthogonal to $h$.

Conversely, if there is a common isotropic vector for all forms in the 2-dimensional subspace $G$, i.e., if all forms in $G$ can be expressed as matrices $\begin{pmatrix} 0 & b \\ b & c \end{pmatrix}$ relative to a suitable basis of $E$, then $G$ also contains the degenerate form represented by the matrix $\begin{pmatrix} 0 & 0 \\ 0 & x \end{pmatrix}$. That form is orthogonal to all forms in *G*, and therefore *G* is not regular under the determinant pairing. ∎

We can now prove the generalization of Desargues' involution theorem stated at the beginning of this note.

**Theorem (Generalization of Desargues' Involution Theorem):** *Consider a projective space of any dimension over a field K of characteristic ≠2. A pencil of quadrics in that projective space will intersect a line in that space (when they do) in pairs of points that are conjugate under an involution if and only if (i) the line is not part of a quadric in that pencil and (ii) the line contains no common zero point for all the quadrics in that pencil. Equivalently, these conditions are satisfied if and only if the pencil contains two quadrics whose intersection points with the given line over an algebraic closure of K are disjoint.*

*Proof.* Consider a pencil of quadrics in the given projective space $P(V)$ defined by a family of symmetric bilinear forms on $V$ of linear dimension 2. Let $\ell$ be a line in that projective space with $\ell = P(E)$ for some 2-dimensional *K*-vector subspace $E$ of $V$. The family of symmetric bilinear forms when restricted to $E$ will give us a linear system of symmetric bilinear forms on $E$ of dimension at most 2.

For the linear system of bilinear forms to be 2-dimensional when restricted to $E$, it is necessary and sufficient that there is no non-zero form in that system that becomes zero when restricted to $E$, i.e., if and only if for any non-zero form $q$ in that system, we do not have $q(u, v) = 0$ for all $u, v$ in $E$, or equivalently, that we do not have $q(v, v) = 0$ for all $v$ in $E$. That is the case if and only if the line $\ell$ is not part of any quadric in the pencil.

Assuming that is the case, the quadrics in the given pencil intersect the line $\ell$ (when they do) in pairs of points conjugate under an involution of the line precisely when the linear system of symmetric bilinear forms as restricted to $E$ is regular under the determinant pairing (in whatever coordinates of the line). In light of Proposition 2, that is equivalent to the lack of a common isotropic point on the line for all the bilinear forms in that linear system, i.e., the line $\ell$ contains no common zero point for all the quadrics in the pencil.

Equivalently, note that the given pencil of quadrics gives us a regular 2-dimensional space of symmetric bilinear forms on $E$, as a vector space over $K$, if and only if that pencil gives us a regular 2-dimensional space of symmetric bilinear forms on $E \otimes L$ as a vector space over an algebraic closure $L$ of $K$. That is because the bilinear forms and their determinant pairing vary naturally with the extension in scalars from $K$ to an algebraic closure of $K$.

If the quadrics give us a regular 2-dimensional space of symmetric bilinear forms on $E \otimes L$, then it is clear that we must have two different quadrics whose intersection points with $P(E \otimes L)$ are disjoint, in light of Proposition 2. Conversely, suppose we have two different quadrics whose intersection points with $P(E \otimes L)$ are disjoint. Note that over an algebraic closure $L$ of $K$ any bilinear form on $E \otimes L$ has isotropic vectors, i.e, any quadric intersects a line in the same projective space. Let $R$ and $S$ are two symmetric bilinear forms that represent these two different quadrics. By the disjoint intersection condition, we can choose two different points $M$ and $N$ of $P(E \otimes L)$, where $M$ is an isotropic point of $R$ but not of $S$ and $N$ is an isotropic point of $S$ but not of $R$. In that case, $M$ and $N$ cannot both be isotropic points for a linear combination $aR + bS$ unless both $a$ and $b$ are zero. In other words, the line $P(E \otimes L)$ cannot be part of any quadric in the given pencil. According to Proposition 2, these two quadrics then give us a regular 2-dimensional space of symmetric bilinear forms on $E \otimes L$. ∎

5.  **GENERAL FORM OF BUTTERFLY THEOREM.** Suppose that $\ell = P(E)$ is a line in a projective space $P(V)$, and we have a Desargues involution on $\ell$ induced by a pencil of quadrics in space $P(V)$. Such an involution will either have two distinct fixed points or no fixed point at all. We have the following characterization of the fixed points of such an involution.

**Proposition 3:** *Let $M$ and $N$ be two distinct points on a line $\ell$ in a projective space $P(V)$, and suppose we have a Desargues involution on $\ell$ induced by a pencil $\mathcal{F}$ of quadrics in $P(V)$. $M$ and $N$ are the fixed points of such a Desargues involution if and only if they are orthogonal relative to all the quadrics in the pencil $\mathcal{F}$.*

*Proof.* Let $\ell = P(E)$ where $E$ is a subspace of dimension 2 of $V$. If $M$ is a fixed point of the Desargues involution on $\ell$ induced by the pencil $\mathcal{F}$, then there is a symmetric bilinear form in that pencil whose restriction to $E$ is degenerate, with $M$ as the radical. Similarly for $N$. The points $M$ and $N$ are orthogonal with respect to these two degenerate bilinear forms, and therefore with respect to all the bilinear forms in the pencil $\mathcal{F}$ because these forms can be expressed as linear combination of any two forms representing two distinct quadrics in the pencil.

Conversely, assume that $M$ and $N$ are orthogonal with respect to all the quadrics in the pencil. Suppose a quadric in the pencil $\mathcal{F}$ intersects $\ell$ in two distinct points $S$ and $T$ (so that $S$ and $T$ are conjugate points under the Desargues involution). Pick vectors $\boldsymbol{m}, \boldsymbol{n}, \boldsymbol{s}, \boldsymbol{t}$ in $E$ representing the points $M, N, S, T$. We can write $\boldsymbol{m} = \alpha\boldsymbol{s} + \beta\boldsymbol{t}$ and $\boldsymbol{n} = \lambda\boldsymbol{s} + \mu\boldsymbol{t}$. If $q$ is a bilinear form representing this quadric, we have $q(\boldsymbol{s}, \boldsymbol{s}) = 0 = q(\boldsymbol{t}, \boldsymbol{t})$ and $q(\boldsymbol{s}, \boldsymbol{t}) \neq 0$ (otherwise the bilinear form $q$ would be zero on $E$). Therefore $q(\boldsymbol{m}, \boldsymbol{n}) = (\alpha\mu + \beta\lambda)q(\boldsymbol{s}, \boldsymbol{t})$. That means $q(\boldsymbol{m}, \boldsymbol{n}) = 0$ if and only if $\alpha\mu + \beta\lambda = 0$, or $\alpha\mu = -\beta\lambda$, i.e., if and only if the cross ratio $(S, T; M, N) = -1$.

Consider the unique involution on $\ell$ with $M$ and $N$ as the fixed points. For any two distinct points $S$ and $T$ that are conjugate under this involution, the cross ratio $(S, T; M, N) = -1$. Accordingly, this involution must be the same as the Desargues involution above, and $M$ and $N$ must therefore be the fixed points of the Desargues involution. ∎

We can now state a general form of the butterfly theorem for an affine space of any dimension over a field $K$ of characteristic $\neq 2$.

**Proposition 4** (General Form of Butterfly Theorem): *Let $\ell$ be a line in an affine space A of any dimension over a field of K of characteristic $\neq 2$. Suppose we have a pencil of quadrics in that affine space and a point M on $\ell$ such that:*

- *A quadric in the pencil is tangent to the line $\ell$ at M, and another quadric in the pencil intersects the line $\ell$ (either over K or an algebraic closure of K) in two distinct points P and Q symmetric about M (i.e., M is the midpoint of P and Q); or*

- *Two different quadrics in the pencil intersect the line ℓ (either over K or an algebraic closure of K) in two different pairs of distinct points, each pair of which is symmetric about M.*

*If either of the above conditions is met, we have the following situation for any quadric in the given pencil: (i) the line ℓ is either tangent to the quadric at M or an asymptote of the quadric (meaning that the line ℓ will intersect and be tangent to the quadric at infinity), or (ii) the line ℓ intersects the quadric in two actual points (over K) or virtual points (over an algebraic closure of K) that are symmetric about M.*

*Proof.* We can regard the affine space $A$ in question as a subset of a projective space $P(V)$ of the same dimension, with the complement of $A$ in $P(V)$ being a hyperplane of $P(V)$, which we regard as the hyperplane at infinity. We can extend ℓ to a projective line $P(E)$ in $P(V)$ by adding a suitable point $N$ to ℓ from the hyperplane at infinity.

The hypothesis in the Proposition implies that the pencil of quadrics induces a Desargues involution on the projective line ℓ ∪ {$N$}, according to our generalization of the Desargues' involution theorem.

Moreover, we claim that the points $M$ and $N$ are orthogonal to each other relative to the two given quadrics, and hence relative to all quadrics in the pencil. To see this, we can work over an algebraic closure of $K$. If a quadric is tangent to the line ℓ at $M$, then $M$ is orthogonal to all the points of ℓ relative to this quadric. If a quadric intersects the line ℓ in two distinct points $P$ and $Q$ symmetric about $M$, then the cross ratio $(P, Q; M, N) = -1$. The proof of Proposition 3 shows that this condition is equivalent to $M$ and $N$ being orthogonal relative to this quadric.

Therefore by Proposition 3, $M$ and $N$ are the fixed points of the Desargues involution. If we choose a coordinate $x$ on the line ℓ ∪ {$N$} so that the point $M$ has coordinate $x = 0$ and the point $N$ is the point at infinity, then the Desargues involution must be the projective transformation $x \mapsto -x$. The results of Proposition 4 follow. ∎

# REFERENCES

1. Scharlau, W. (1985).  *Quadratic and Hermitian Forms*.  Grundlehren der Mathematischen Wissenschaften, vol. 270.  Berlin, Heidelberg, New York and Tokyo: Springer-Verlag.
2. Pedoe, D. (1988).  *Geometry, A Comprehensive Course*.  New York: Dover.

NICHOLAS PHAT NGUYEN

*12015 12th Dr SE, Everett, WA 98208, U.S.A.*

*Email:  nicholas.pn@gmail.com*